# MINIMUM DISTANCE REGRESSION MODEL CHECKING WITH BERKSON MEASUREMENT ERRORS

By Hira L. Koul[1] and Weixing Song

*Michigan State University and Kansas State University*

Lack-of-fit testing of a regression model with Berkson measurement error has not been discussed in the literature to date. To fill this void, we propose a class of tests based on minimized integrated square distances between a nonparametric regression function estimator and the parametric model being fitted. We prove asymptotic normality of these test statistics under the null hypothesis and that of the corresponding minimum distance estimators under minimal conditions on the model being fitted. We also prove consistency of the proposed tests against a class of fixed alternatives and obtain their asymptotic power against a class of local alternatives orthogonal to the null hypothesis. These latter results are new even when there is no measurement error. A simulation that is included shows very desirable finite sample behavior of the proposed inference procedures.

**1. Introduction.** A classical problem in statistics is to use a vector $X$ of $d$-dimensional variables, $d \geq 1$, to explain the one-dimensional response $Y$. As is the practice, this is often done in terms of the regression function $\mu(x) := E(Y|X=x), x \in \mathbb{R}^d$, assuming it exists. Usually, in practice the predictor vector $X$ is assumed to be observable. But in many experiments, it is expensive or impossible to observe $X$. Instead, a proxy or a manifest $Z$ of $X$ can be measured. As an example, consider the herbicide study of Rudemo, Ruppert and Streibig [16] in which a nominal measured amount $Z$ of herbicide was applied to a plant but the actual amount absorbed $X$ by the plant is unobservable. As another example, from Wang [20], an epidemiologist studies the severity of a lung disease, $Y$, among the residents in a city in relation to the amount of certain air pollutants, $X$. The amount of the air pollutants $Z$ can be measured at certain observation stations in the city, but the actual exposure of the residents to the pollutants, $X$, is

Received October 2006.
[1]Supported in part by NSF Grant DMS-07-04130.

unobservable and may vary randomly from the $Z$-values. In both cases, $X$ can be expressed as $Z$ plus a random error. There are many similar examples in agricultural or medical studies; see, for example, Carroll, Ruppert and Stefanski [5] and Fuller [10], among others. All these examples can be formalized into the so-called Berkson model

$$(1.1) \qquad Y = \mu(X) + \varepsilon, \qquad X = Z + \eta,$$

where $\eta$ and $\varepsilon$ are random errors with $E\varepsilon = 0$, $\eta$ is $d$-dimensional and $Z$ is the observable $d$-dimensional control variable. All three r.v.'s $\varepsilon, \eta$ and $Z$ are assumed to be mutually independent.

Let $\mathcal{M} := \{m_\theta(x) : x \in \mathbb{R}^d, \theta \in \Theta \subset \mathbb{R}^q\}$, $q \geq 1$, be a class of known functions. The parametric Berkson regression model where $\mu \in \mathcal{M}$ has been the focus of numerous authors. Cheng and Van Ness [6] and Fuller [10], among others, discuss the estimation in the linear Berkson model. For nonlinear models, [5] and references therein consider the estimation problem by using a regression calibration method. Huwang and Huang [13] study the estimation problem when $m_\theta(x)$ is a polynomial in $x$ of a known order and show that the least square estimators based on the first two conditional moments of $Y$, given $Z$, are consistent. Similar results are obtained in [19] and [20] for a class of nonlinear Berkson models.

But literature appears to be scant on the lack-of-fit testing problem in this important model. This paper makes an attempt in filling this void. To be precise, with $(Z, Y)$ obeying the model (1.1), the problem of interest here is to test the hypothesis

$\mathcal{H}_0 : \mu(x) = m_{\theta_0}(x) \qquad$ for some $\theta_0 \in \Theta$ and for all $x \in \mathcal{I}$;

$\mathcal{H}_1 : \mathcal{H}_0$ is not true,

based on a random sample $(Z_i, Y_i), 1 \leq i \leq n$, from the distribution of $(Z, Y)$, where $\Theta$ and $\mathcal{I}$ are compact subsets of $\mathbb{R}^q$ and $\mathbb{R}^d$, respectively.

Interesting and profound results, on the contrary, are available for regression model checking in the absence of errors in independent variables; see, for example, [1, 11, 12] and references therein, [17, 18], among others. Koul and Ni [14] use the minimum distance methodology to propose tests of lack-of-fit of a parametric regression model in the classical regression setup. In a finite sample comparison of these tests with some other existing tests, they noted that a member of this class preserves the asymptotic level and has relatively very high power against some alternatives. The present paper extends this methodology to the above Berkson model.

To be specific, Koul and Ni considered the following tests of $\mathcal{H}_0$ where the design is random and observable, and the errors are heteroscedastic.

For any density kernel $K$, let $K_h(x) := K(x/h)/h^d$, $h > 0, x \in \mathbb{R}^d$. Define $\tilde{f}_w(x) := \frac{1}{n} \sum_{j=1}^n K_w^*(x - X_j), w = w_n \sim (\log n/n)^{1/(d+4)}$,

$$T_n(\theta) := \int_{\mathcal{C}} \left[ \frac{1}{n} \sum_{j=1}^n K_h(x - X_j)(Y_j - m_\theta(X_j)) \right]^2 \frac{d\bar{G}(x)}{\tilde{f}_w^2(x)}$$

and $\tilde{\theta}_n := \arg\min\{T_n(\theta), \theta \in \Theta\}$, where $K, K^*$ are density kernel functions, possibly different, $h = h_n$ and $w = w_n$ are the window widths, depending on the sample size $n$, $\mathcal{C}$ is a compact subset of $\mathbb{R}^d$ and $\bar{G}$ is a $\sigma$-finite measure on $\mathcal{C}$. They proved consistency and asymptotic normality of this estimator, and that the asymptotic null distribution of $\mathcal{D}_n := nh_n^{d/2}(T_n(\tilde{\theta}_n) - \tilde{C}_n)/\tilde{\Gamma}_n^{1/2}$ is standard normal, where

$$\tilde{C}_n := n^{-2} \sum_{i=1}^n \int_{\mathcal{C}} K_h^2(x - X_i)\varepsilon_i^2 \tilde{f}_w^{-2}(x) \, d\bar{G}(x), \qquad \hat{\varepsilon}_i = Y_i - m_{\tilde{\theta}_n}(X_i),$$

$$\tilde{\Gamma}_n := h^d n^{-2} \sum_{i \neq j=1}^n \left( \int_{\mathcal{C}} K_h(x - X_i)K_h(x - X_j)\hat{\varepsilon}_i\hat{\varepsilon}_j \tilde{f}_w^{-2}(x) \, d\bar{G}(x) \right)^2.$$

These results were made feasible by recognizing to use an optimal window width $w_n$ for the estimation of the denominator and a different window width $h_n$ for the estimation of the numerator in the kernel-type nonparametric estimator of the regression function. A consequence of the above asymptotic normality result is that at least for large samples one does not need to use any resampling method to implement these tests.

These findings thus motivate one to look for tests of lack-of-fit in the Berkson model based on the above minimized distances. Since the predictors in Berkson models are unobservable, clearly the above procedures need some modifications.

Let $f_\varepsilon$, $f_X$, $f_\eta$, $f_Z$ denote the density functions of the r.v.'s in their subscripts and $\sigma_\varepsilon^2$ denote the variance of $\varepsilon$. In linear regression models if one is interested in making inference about the coefficient parameters only, these density functions need not be known. Berkson [3] pointed out that the ordinary least square estimators are unbiased and consistent in these models and one can simply ignore the measurement error $\eta$. But if the regression model is nonlinear or if there are other parameters in the Berkson model that need to be estimated, then extra information about these densities should be supplied to ensure the identifiability. A standard assumption in the literature is to assume that $f_\eta$ is known or unknown up to a Euclidean parameter; compare [5, 13, 20], among others. For the sake of relative transparency of the exposition we assume that $f_\eta$ is known.

To adopt the Koul and Ni (K–N) procedure to the current setup, we first need to obtain a nonparametric estimator of $\mu$. Note that in the model (1.1), $f_X(x) = \int f_Z(z) f_\eta(x-z)\, dz$. For any kernel density $K$, let $K_{hi}(z) := K_h(z - Z_i)$, $\hat{f}_{Zh}(z) = \sum_{i=1}^n K_{hi}(z)/n$ and $\bar{K}_h(x, z) := \int K_h(z-y) f_\eta(x-y)\, dy$, for $x, z \in \mathbb{R}^d$. It is then natural to estimate $f_X(x)$ and $\mu(x)$ by

$$\hat{f}_X(x) := \frac{1}{n} \sum_{i=1}^n \bar{K}_h(x, Z_i), \qquad \hat{J}_n(x) := \frac{\sum_{i=1}^n \bar{K}_h(x, Z_i) Y_i}{\sum_{i=1}^n \bar{K}_h(x, Z_i)}.$$

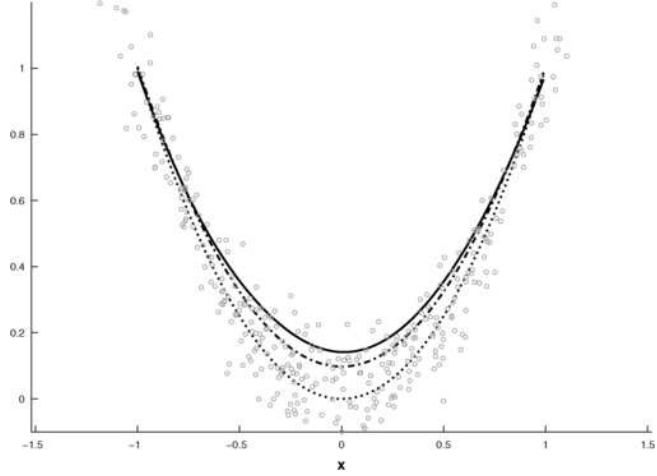

Fig. 1.

A routine argument, however, shows that $\hat{J}_n(x)$ is a consistent estimator of $J(x) := E[H(Z)|X = x]$, where $H(z) := E[\mu(X)|Z = z]$, but not of $\mu(x)$.

We include the following simulation study to illustrate this point. Consider the model $Y = X^2 + \varepsilon$, $X = Z + \eta$, where $\varepsilon$ and $\eta$ are Gaussian r.v.'s with means zero and variances 0.01 and 0.05, respectively, and $Z$ is a standard Gaussian r.v. Then $J(x) = 0.0976 + 0.907 x^2$. We generated 500 samples from this model, calculated $\hat{J}_n$ and then put all three graphs, $\hat{J}_n(x)$, $\mu(x) = x^2$, $J(x) = 0.0976 + 0.907 x^2$ into one plot in Figure 1. The curves with solid, dash-dot and dot lines are those of $\hat{J}_n$, $J(x)$ and $\mu(x) = x^2$, respectively.

To overcome this difficulty, one way to proceed is as follows. Define

$$H_\theta(z) := E[m_\theta(X)|Z = z], \qquad J_\theta(x) := E[H_\theta(Z)|X = x],$$

$$(1.2) \qquad \widetilde{Q}_n(\theta) = \int_{\mathcal{C}} \left[ \frac{1}{n \hat{f}_X(x)} \sum_{i=1}^n \bar{K}_h(x, Z_i) Y_i - J_\theta(x) \right]^2 d\bar{G}(x),$$

$$Q_n(\theta) = \int_{\mathcal{C}} \left[ \frac{1}{n \hat{f}_X(x)} \sum_{i=1}^n \bar{K}_h(x, Z_i)[Y_i - H_\theta(Z_i)] \right]^2 d\bar{G}(x)$$

and $\widetilde{\theta}_n = \arg\min_{\theta \in \Theta} \widetilde{Q}_n(\theta), \theta_n = \arg\min_{\theta \in \Theta} Q_n(\theta)$.

Under some conditions, we can show that $\theta_n$, $\widetilde{\theta}_n$ are consistent for $\theta$ and asymptotic null distribution of a suitably standardized $Q_n(\theta_n)$ is the same as that of a degenerate $U$-statistic, whose asymptotic distribution in turn is the same as that of an infinite sum of weighted centered chi-square random variables. Since the kernel function in the degenerate $U$-statistic is complicated, computation of its eigenvalues and eigenfunctions is not easy and hence this test is hard to implement in practice.

An alternative way to proceed is to use regression calibration as follows. Because $E(Y|Z) = H(Z)$, one considers the new regression model $Y = H(Z) + \zeta$, where the error $\zeta$ satisfies $E(\zeta|Z) = 0$. The problem of testing for $\mathcal{H}_0$ is now transformed to that of testing for $H(z) = H_{\theta_0}(z)$. This motivates the following modification of the K–N procedure that adjusts for not observing the design variable. Let

$$\hat{f}_{Zw}(z) := \frac{1}{n}\sum_{i=1}^{n} K^*_{wi}(z), \qquad \hat{H}_n(z) := \frac{\sum_{i=1}^{n} K_{hi}(z)Y_i}{n\hat{f}_{Zw}(z)}, \qquad z \in \mathbb{R}^d.$$

Note that $\hat{H}_n$ is an estimator of $H(z) = E(\mu(X)|Z=z)$. Define

$$M_n^*(\theta) = \int_{\mathcal{I}} \left[\frac{1}{n\hat{f}_{Zw}(z)} \sum_{i=1}^{n} K_{hi}(z)Y_i - H_\theta(z)\right]^2 dG(z),$$

(1.3) $$M_n(\theta) = \int_{\mathcal{I}} \left[\frac{1}{n\hat{f}_{Zw}(z)} \sum_{i=1}^{n} K_{hi}(z)[Y_i - H_\theta(Z_i)]\right]^2 dG(z),$$

$$\theta_n^* = \arg\min_{\theta \in \Theta} M_n^*(\theta), \qquad \hat{\theta}_n = \arg\min_{\theta \in \Theta} M_n(\theta),$$

where $G$ is a $\sigma$-finite measure supported on $\mathcal{I}$. We consider $M_n$ to be the right analog of the above $T_n$ for the Berkson model. This paper establishes consistency of $\theta_n^*$ and $\hat{\theta}_n$ for $\theta_0$ and asymptotic normality of $\sqrt{n}(\hat{\theta}_n - \theta_0)$, under $\mathcal{H}_0$. Additionally, we prove that the asymptotic null distribution of the normalized test statistic $\widehat{\mathcal{D}}_n := nh^{d/2}\hat{\Gamma}_n^{-1/2}(M_n(\hat{\theta}_n) - \hat{C}_n)$ is standard normal, which, unlike the modification (1.2), can be easily used to implement this testing procedure, at least for the large samples. Here,

$$d\hat{\psi}(z) := \frac{dG(z)}{\hat{f}_{Zw}^2(z)}, \qquad \hat{\zeta}_i := Y_i - H_{\hat{\theta}_n}(Z_i), \qquad 1 \leq i \leq n,$$

(1.4) $$\hat{C}_n := \frac{1}{n^2}\sum_{i=1}^{n} \int K_{hi}^2(z)\hat{\zeta}_i^2\, d\hat{\psi}(z),$$

$$\hat{\Gamma}_n := \frac{2h^d}{n^2}\sum_{i\neq j}\left(\int K_{hi}(z)K_{hj}(z)\hat{\zeta}_i\hat{\zeta}_j\, d\hat{\psi}(z)\right)^2.$$

We note that a factor of 2 is missing in the analog of $\hat{\Gamma}_n$ in K–N.

Even though K–N conducted some convincing simulations to demonstrate the finite sample power properties of the $\mathcal{D}_n$-tests, they did not discuss any theoretical power properties of their tests. In contrast, we prove consistency of the proposed minimum distance (MD) tests against a large class of fixed alternatives and obtain their asymptotic power under a class of local alternatives. Let $L_2(G)$ denote the class of real-valued square integrable functions on $\mathbb{R}^d$ with respect to $G$, $\rho(\nu_1, \nu_2) := \int [\nu_1 - \nu_2]^2 \, dG, \nu_1, \nu_2 \in L_2(G)$ and

(1.5) $$T(\nu) := \arg\min_{\theta \in \Theta} \rho(\nu, H_\theta), \qquad \nu \in L_2(G).$$

Let $m \in L_2(G)$ be a given function. Consider the problem of testing $\mathcal{H}_0$ against the alternative $\mathcal{H}_1 : \mu = m, m \notin \mathcal{M}$. Under assumption (m2) below and $\mathcal{H}_0$, $T(H_{\theta_0}) = \theta_0$, while under $\mathcal{H}_1$, $T(H) \neq \theta_0$, where now $H(z) = E(m(X)|Z = z)$. Consistency of the $\widehat{\mathcal{D}}_n$-test requires consistency of $\hat{\theta}_n$ for $T(H)$ only, while its asymptotic power properties against the local alternatives $\mathcal{H}_{1n} : \mu = m_{\theta_0} + r/nh^{d/2}$ requires that $n^{1/2}(\hat{\theta}_n - \theta_0) = O_p(1)$, under $\mathcal{H}_{1n}$. Here $r$ is a continuously differentiable function with $R(z) := E(r(X)|Z = z)$ such that $R \in L_2(G)$ and $\int H_\theta R \, dG = 0$ for all $\theta \in \Theta$. Under assumptions of Section 2 below, we show that under $\mathcal{H}_1$, $\hat{\theta}_n \to T(H)$, in probability, and under $\mathcal{H}_{1n}$, both $n^{1/2}(\hat{\theta}_n - \theta_0)$ and $\widehat{\mathcal{D}}_n$ are asymptotically normally distributed.

The paper is organized as follows. The needed assumptions are stated in the next section. All limits are taken as $n \to \infty$, unless mentioned otherwise. Section 3 contains the proofs of consistency of $\theta_n^*$ and $\hat{\theta}_n$ while Section 4 discusses asymptotic normality of $\hat{\theta}_n$ and $\widehat{\mathcal{D}}_n$, under $\mathcal{H}_0$. The power of the MD-test for fixed and local alternatives is discussed in Section 5. The simulation results in Section 6 show little bias in the estimator $\hat{\theta}_n$ for all chosen sample sizes. The finite sample level approximates the nominal level well for larger sample sizes and the empirical power is high (above 0.9) for moderate to large sample sizes against the chosen alternatives.

Finally, we mention that closely related to the Berkson model is the so-called errors-in-variable regression model in which $Z = X + u$. In this case also one can use the above MD method to test $\mathcal{H}_0$, although we do not deal with this here. The biggest challenge is to construct nonparametric estimators of $f_X$ and $H_\theta$. The deconvolution estimators discussed in Fan [7, 8], Fan and Truong [9], among others, may be found useful here.

**2. Assumptions.** Here we shall state the needed assumptions in this paper. In the assumptions below $\theta_0$ denotes the true parameter value under $\mathcal{H}_0$. About the errors, the underlying design and $G$ we assume the following:

(e1) $\{(Z_i, Y_i) : Z_i \in \mathbb{R}^d, i = 1, 2, \ldots, n\}$ are i.i.d. with $H(z) := E(Y|Z = z)$ satisfying $\int H^2 \, dG < \infty$, where $G$ is a $\sigma$-finite measure on $\mathcal{I}$.

(e2) $0 < \sigma_\varepsilon^2 < \infty$, $Em_{\theta_0}^2(X) < \infty$ and the function $\tau^2(z) = E[(m_{\theta_0}(X) - H_{\theta_0}(Z))^2|Z=z]$ is a.s. (G) continuous on $\mathcal{I}$.
(e3) $E|\varepsilon|^{2+\delta} < \infty$, $E|m_{\theta_0}(X) - H_{\theta_0}(Z)|^{2+\delta} < \infty$, for some $\delta > 0$.
(e4) $E|\varepsilon|^4 < \infty$, $E|m_{\theta_0}(X) - H_{\theta_0}(Z)|^4 < \infty$.
(f1) $f_Z$ is uniformly continuous and bounded from below on $\mathcal{I}$.
(f2) $f_Z$ is twice continuously differentiable.
(g) $G$ has a continuous Lebesgue density $g$ on $\mathcal{I}$.

About the bandwidth $h_n$ we shall make the following assumptions:

(h1) $h_n \to 0$.
(h2) $nh_n^{2d} \to \infty$.
(h3) $h_n \sim n^{-a}$, where $0 < a < \min(1/2d, 4/(d(d+4)))$.

About the kernel functions $K$ and $K^*$ we shall assume the following:

(k) The kernel functions $K$, $K^*$ are positive symmetric square integrable densities on $[-1, 1]^d$. In addition, $K^*$ satisfies a Lipschitz condition.

About the parametric family $\{m_\theta\}$ we assume the following:

(m1) For each $\theta$, $m_\theta(x)$ is a.e. continuous in $x$ w.r.t. the Lebesgue measure.
(m2) The function $H_\theta(z)$ is identifiable w.r.t. $\theta$, that is, if $H_{\theta_1}(z) = H_{\theta_2}(z)$ for almost all $z(G)$, then $\theta_1 = \theta_2$.
(m3) For some positive continuous function $\ell$ on $\mathcal{I}$ and for some $0 < \beta \leq 1$, $|H_{\theta_2}(z) - H_{\theta_1}(z)| \leq \|\theta_2 - \theta_1\|^\beta \ell(z)$, $\forall \theta_1, \theta_2 \in \Theta, z \in \mathcal{I}$.

For every $z$, $H_\theta(z)$ is differentiable in $\theta$ in a neighborhood of $\theta_0$ with the vector of derivative $\dot{H}_\theta(z)$ satisfying the following three conditions:

(m4) $\forall 0 < \delta_n \to 0$

$$\sup_{1 \leq i \leq n, \|\theta - \theta_0\| \leq \delta_n} \frac{|H_\theta(Z_i) - H_{\theta_0}(Z_i) - (\theta - \theta_0)'\dot{H}_{\theta_0}(Z_i)|}{\|\theta - \theta_0\|} = o_p(1).$$

(m5) $\forall 0 < k < \infty$

$$\sup_{1 \leq i \leq n, \sqrt{nh_n^d}\|\theta - \theta_0\| \leq k} h_n^{-d/2} \|\dot{H}_\theta(Z_i) - \dot{H}_{\theta_0}(Z_i)\| = o_p(1).$$

(m6) $\int \|\dot{H}_{\theta_0}\|^2 \, dG < \infty$ and $\Sigma_0 := \int \dot{H}_{\theta_0} \dot{H}'_{\theta_0} \, dG$ is positive definite.

For later use we note that, under (h2) and (m4), $nh^d \to \infty$ and for every $0 < k < \infty$,

$$(2.1) \quad \sup_{1 \leq i \leq n, \sqrt{nh_n^d}\|\theta - \theta_0\| \leq k} \frac{|H_\theta(Z_i) - H_{\theta_0}(Z_i) - (\theta - \theta_0)'\dot{H}_{\theta_0}(Z_i)|}{\|\theta - \theta_0\|} = o_p(1).$$

The above conditions are similar to those imposed in K–N on the model $m_\theta$. Consider the following conditions in terms of the given model:

(m2′) The parametric family of models $m_\theta(x)$ is identifiable w.r.t. $\theta$, that is, if $m_{\theta_1}(x) = m_{\theta_2}(x)$ for almost all $x$, then $\theta_1 = \theta_2$.

(m3′) For some positive continuous function $L$ on $\mathbb{R}^d$ with $EL(X) < \infty$ and for some $\beta > 0$, $|m_{\theta_2}(x) - m_{\theta_1}(x)| \leq \|\theta_2 - \theta_1\|^\beta L(x)$, $\forall \theta_1, \theta_2 \in \Theta, x \in \mathbb{R}^d$.

The function $m_\theta(x)$ is differentiable in $\theta$ in a neighborhood of $\theta_0$, with the vector of differential $\dot{m}_{\theta_0}$ satisfying the following two conditions:

(m4′) $\forall 0 < \delta_n \to 0$

$$\sup_{x \in \mathbb{R}^d, \|\theta - \theta_0\| \leq \delta_n} \frac{|m_\theta(x) - m_{\theta_0}(x) - (\theta - \theta_0)' \dot{m}_{\theta_0}(x)|}{\|\theta - \theta_0\|} \to 0.$$

(m5′) For every $0 < k < \infty$

$$\sup_{x \in \mathbb{R}^d, \sqrt{nh_n^d}\|\theta - \theta_0\| \leq k} h_n^{-d/2} \|\dot{m}_\theta(x) - \dot{m}_{\theta_0}(x)\| = o(1).$$

In some cases, (m2) and (m2′) are equivalent. For example, if the family of densities $\{f_\eta(\cdot - z); z \in \mathbb{R}\}$ is complete, then this holds. Similarly, if $m_\theta(x) = \theta'\gamma(x)$ and $\int \gamma(x) f_\eta(x-z)\,dx \neq 0$, for all $z$, then also (m2) and (m2′) are equivalent.

We can also show that (m3′)–(m5′) imply (m3)–(m5), respectively. This follows because $H_\theta(z) \equiv \int m_\theta(x) f_\eta(x-z)\,dx$, so that under (m3′), $|H_{\theta_2}(z) - H_{\theta_1}(z)| \leq \|\theta_2 - \theta_1\|^\beta \int L(x) f_\eta(x-z)\,dx, \forall z \in \mathbb{R}^d$. Hence (m3) holds with $\ell(z) = \int L(x) f_\eta(x-z)\,dx$. Note that $E\ell(Z) = EL(X) < \infty$.

Using the fact that $\int f_\eta(x-z)\,dx \equiv 1$, the left-hand side of (m4) is bounded above by $\sup_{x \in \mathbb{R}^d, \|\theta - \theta_0\| \leq \delta} |m_\theta(x) - m_{\theta_0}(x) - (\theta - \theta_0)' \dot{m}_{\theta_0}(x)|/\|\theta - \theta_0\| = o(1)$, by (m4′). Similarly, (m5′) implies (m5) and (m1) implies that $H_\theta(z)$ is a.s. continuous in $z(G)$.

The conditions (m1)–(m6) are trivially satisfied when $m_\theta(x) = \theta'\gamma(x)$ provided components of $E[\gamma(X)|Z = z]$ are continuous, nonzero on $\mathcal{I}$ and the matrix $\int E[\gamma(X)\gamma'(X)|Z = z]\,dG(z)$ is positive definite.

The conditions (e1), (e2), (f1), (k), (m1)–(m3), (h1) and (h2) suffice for consistency of $\hat{\theta}_n$, while these plus (e3), (f2), (m4)–(m6) and (h3) are needed for the asymptotic normality of $\hat{\theta}_n$. The asymptotic normality of $M_n(\hat{\theta}_n)$ needs (e1)–(e4) and (f1)–(m6) and (h3). Of course, (h3) implies (h1) and (h2).

Let $q_{h1} := f_Z/\hat{f}_{Zh} - 1$. From [15] we obtain that under (f1), (k), (h1) and (h2),

$$\sup_{z \in \mathcal{I}} |\hat{f}_{Zh}(z) - f_Z(z)| = o_p(1), \qquad \sup_{z \in \mathcal{I}} |\hat{f}_{Zw}(z) - f_Z(z)| = o_p(1),$$

(2.2)
$$\sup_{z \in \mathcal{I}} |q_{h1}(z)| = o_p(1) = \sup_{z \in \mathcal{I}} |q_{w1}(z)|.$$

These conclusions are often used in the proofs below.

In the sequel, the true parameter $\theta_0$ is assumed to be an inner point of $\Theta$ and $\zeta := Y - H_{\theta_0}(Z)$. The integrals with respect to $G$ are understood to be over $\mathcal{I}$. The convergence in distribution is denoted by $\to_d$ and $\mathcal{N}_p(a, B)$ denotes the $p$-dimensional normal distribution with mean vector $a$ and covariance matrix $B$, $p \geq 1$. We shall also need the following notation:

$$d\psi(z) := \frac{dG(z)}{f_Z^2(z)}, \qquad \sigma_\zeta^2(z) := \operatorname{Var}_{\theta_0}(\zeta|Z=z) = \sigma_\varepsilon^2 + \tau^2(z),$$

$$\zeta_i := Y_i - H_{\theta_0}(Z_i), \qquad 1 \leq i \leq n,$$

$$\tilde{C}_n := n^{-2} \sum_{i=1}^n \int K_{hi}^2 \zeta_i^2 \, d\psi,$$

(2.3)
$$K_2(v) := \int K(v+u)K(u) \, du, \qquad \|K_2\|^2 := \int K_2^2(v) \, dv,$$

$$\Gamma := 2\|K_2\|^2 \int (\sigma_\zeta^2(z))^2 g(z) \, d\psi(z),$$

$$q_n(z) := (f_Z^2(z)/\hat{f}_{Zw}^2(z)) - 1,$$

$$\mu_n(z, \theta) := \frac{1}{n} \sum_{i=1}^n K_{hi}(z) H_\theta(Z_i), \qquad \dot{\mu}_n(z, \theta) := \frac{1}{n} \sum_{i=1}^n K_{hi}(z) \dot{H}_\theta(Z_i),$$

(2.4) $\quad U_n(z, \theta) := \dfrac{1}{n} \sum_{i=1}^n K_{hi}(z)[Y_i - H_\theta(Z_i)], \qquad U_n(z) := U_n(z, \theta_0),$

$$\mathcal{Z}_n(z, \theta) := \frac{1}{n} \sum_{i=1}^n K_{hi}(z)[H_\theta(Z_i) - H_{\theta_0}(Z_i)], \qquad \theta \in \mathbb{R}^q, z \in \mathbb{R}^d.$$

These entities are analogous to the similar entities defined at (3.1) in K–N. The main difference is that $\mu_\theta$ there is replaced by $H_\theta$ and $X_i$'s by $Z_i$'s.

**3. Consistency of $\theta_n^*$ and $\hat{\theta}_n$.** Recall (1.5). In this section we first prove consistency of $\theta_n^*$ and $\hat{\theta}_n$ for $T(H)$, where $H$ corresponds to a given regression function $m$. Consistency of these estimators for $\theta_0$ under $\mathcal{H}_0$ follows from this general result. The following lemma is found useful in the proofs here. Its proof is similar to that of Theorem 1 in [2].

LEMMA 3.1. *Under the conditions* (m3), *the following hold:*
(a) $T(\nu)$ *always exists, for all* $\nu \in L_2(G)$.
(b) *If* $T(\nu)$ *is unique, then $T$ is continuous at $\nu$ in the sense that for any sequence of* $\{\nu_n\} \in L_2(G)$ *converging to $\nu$ in $L_2(G)$, $T(\nu_n) \to T(\nu)$, that is, $\rho(\nu_n, \nu) \to 0$ implies $T(\nu_n) \to T(\nu)$.*

(c) *In addition, if* (m2) *holds, then* $T(H_\theta) = \theta$, *uniquely for all* $\theta \in \Theta$.

From now on, we use the convention that for any integral $J := \int r\, d\hat{\psi}$, $\tilde{J} := \int r\, d\psi$. Also, let $\gamma^2(z) := E[(m(X) - H(Z))^2 | Z = z]$, $z \in \mathbb{R}^d$. A consequence of the above lemma is the following.

LEMMA 3.2. *Suppose* (k), (f1), (m3) *hold and $m$ is a given regression function satisfying the model assumption* (1.1), $H \in L_2(G)$ *and $T(H)$ is unique.*

(a) *In addition, suppose $H$ and $\gamma^2$ are a.e. (G) continuous. Then,* $\theta_n^* = T(H) + o_p(1)$.

(b) *In addition, suppose $m$ is continuous on $\mathcal{I}$. Then,* $\hat{\theta}_n = T(H) + o_p(1)$.

PROOF.

PROOF OF PART (a). We shall use part (b) of Lemma 3.1 with $\nu_n = \hat{H}_n$ and $\nu = H$. Note that $M_n^*(\theta) = \rho(\hat{H}_n, H_\theta)$, $\theta_n^* = T(\hat{H}_n)$. It thus suffices to prove

$$(3.1) \qquad \rho(\hat{H}_n, H) = o_p(1).$$

Let $\xi_i := Y_i - H(Z_i)$, $1 \le i \le n$,

$$U_n(z) := n^{-1} \sum_{i=1}^n K_{hi}(z) \xi_i,$$

$$(3.2) \qquad \bar{H}(z) := n^{-1} \sum_{i=1}^n K_{hi}(z) H(Z_i), \qquad z \in \mathbb{R}^d,$$

$$\Delta_n := \int [\bar{H} - \hat{f}_{Zw} H]^2 \, d\hat{\psi}.$$

To prove (3.1), plug $Y_i = \xi_i + H(Z_i)$ in $\rho(\hat{H}_n, H)$ and expand the quadratic integrand to obtain that $\rho(\hat{H}_n, H) \le 2[\int U_n^2 \, d\hat{\psi} + \Delta_n]$. By Fubini's theorem and orthogonality of $Z_i$ and $\xi_i$,

$$(3.3) \quad E \int U_n^2(z)\, d\psi(z) = n^{-1} \int E\{K_h^2(z - Z)(\sigma_\varepsilon^2 + \gamma^2(Z))\}\, d\psi(z).$$

By the continuity of $f_Z$ [cf. (f1)], by a.e. continuity of $\gamma^2$ and by (k), we obtain, for $j = 0, 2$, that

$$EK_h^2(z - Z)\gamma^j(Z) = \frac{1}{h^d} \int K^2(y) f_Z(z - yh) \gamma^j(z - yh)\, dy = O\left(\frac{1}{h^d}\right).$$

These calculations, the bound $\int U_n^2 \, d\hat{\psi} \le \sup_{z \in \mathcal{I}} (\frac{f_Z(z)}{\hat{f}_{Zw}(z)})^2 \int U_n^2\, d\psi$ and (2.2) imply that

$$(3.4) \qquad E\int U_n^2\, d\psi = O\left(\frac{1}{nh^d}\right) \quad \text{and} \quad \int U_n^2\, d\hat{\psi} = O_p\left(\frac{1}{nh^d}\right).$$

Next, we shall show that

(3.5) $$\Delta_n = o_p(1).$$

Toward this goal, add and subtract $H(z)E(\hat{f}_{Zw}(z)) = H(z)E(K_h^*(z-Z))$ and $E(\bar{H}(z)) = E(K_h(z-Z)H(Z))$ in the quadratic term of the integrand in $\Delta_n$, to obtain $\Delta_n \leq 4[\Delta_{n1} + \Delta_{n2} + \Delta_{n3}]$, where $\Delta_{n1} := \int [\bar{H} - E(\bar{H})]^2 \, d\hat{\psi}$, $\Delta_{n2} := \int [\hat{f}_{Zw} - E(\hat{f}_{Zw})]^2 H^2 \, d\hat{\psi}$, $\Delta_{n3} := \int [E(\bar{H}) - HE(\hat{f}_{Zw})]^2 \, d\hat{\psi}$.

Fubini's theorem, (k), (f1) and $H$ being a.e. $(G)$ continuous imply

$$E\tilde{\Delta}_{n1} \leq n^{-1} \int E[K_h^2(z-Z)H^2(Z)] \, d\psi(z)$$
$$= (nh^d)^{-1} \int\int K^2(w)H^2(z-wh)f_Z(z-wh) \, dw \, d\psi(z)$$
$$= O((nh^d)^{-1}).$$

Because $\Delta_{n1} \leq \sup_z (f_Z(z)/\hat{f}_{Zw}(z))^2 \tilde{\Delta}_{n1}$, the above bound and (2.2) yield that $\Delta_{n1} = O_p((nh^d)^{-1})$. Similarly, one shows that $\Delta_{n2} = O_p((nh^d)^{-1})$.

Next, $H$ being a.e. $(G)$ continuous and (f1) yield

$$\tilde{\Delta}_{n3} = \int \left[\int [K(u)H(z-hu) - H(z)K^*(u)]f_Z(z-hu) \, du\right]^2 d\psi(z) \to 0.$$

Hence, by (2.2), $\Delta_{n3} = o_p(1)$. This completes the proof of (3.5) and hence that of part (a).

PROOF OF PART (b). Consistency of $\hat{\theta}_n$ for $\theta_0$ under $\mathcal{H}_0$ can be proved by using the method in [14]. But that method does not yield consistency of $\hat{\theta}_n$ for $T(H)$ when $\mu = m$, $m \notin \mathcal{M}$. The proof in general consists of showing

(3.6) $$\sup_{\theta \in \Theta} |M_n(\theta) - \rho(H, H_\theta)| = o_p(1).$$

This, (m3) and the continuity of $m$ on $\mathcal{I}$ imply that $H$ is continuous and $|\rho(H, H_{\theta_2}) - \rho(H, H_{\theta_1})| \leq C\|\theta_1 - \theta_2\|^\beta$, $\forall \theta_1, \theta_2 \in \Theta$, which in turn implies that for all $\epsilon > 0$,

(3.7) $$\lim_{\delta \to 0} \limsup_n P\left(\sup_{\|\theta_1 - \theta_2\| < \delta} |M_n(\theta_1) - M_n(\theta_2)| > \epsilon\right) = 0.$$

These two facts in turn imply $\hat{\theta}_n = T(H) + o_p(1)$. For, suppose $\hat{\theta}_n \not\to T(H)$, in probability. Then, by the compactness of $\Theta$, there is a subsequence $\{\hat{\theta}_{n_k}\}$ of $\{\hat{\theta}_n\}$ and a $\theta^* \neq T(H)$ such that $\hat{\theta}_{n_k} = \theta^* + o_p(1)$. Because $M_{n_k}(\hat{\theta}_{n_k}) \leq M_{n_k}(T(H))$, we obtain

$$\rho(H, H_{\theta^*}) \leq \rho(H, H_{T(H)}) + 2\sup_\theta |M_{n_k}(\theta) - \rho(H, H_\theta)|$$
$$+ |M_{n_k}(\theta^*) - M_{n_k}(\hat{\theta}_{n_k})|.$$

By (3.6) and (3.7), the last two summands in the above bound are $o_p(1)$, so that $\rho(H, H_{\theta^*}) \leq \rho(H, H_{T(H)})$ eventually, with arbitrarily large probability. In view of the uniqueness of $T(H)$, this is a contradiction unless $\theta^* = T(H)$.

To prove (3.6), use the Cauchy–Schwarz (C–S) inequality to obtain that $|M_n(\theta) - \rho(H, H_\theta)|$ is bounded above by the product $Q_{n1}(\theta) Q_{n2}(\theta)$, where

$$Q_{n1}(\theta) := \int \left( [\hat{H}_n(z) - H(z)] - \left[ \frac{\mu_n(z, \theta)}{\hat{f}_{Zw}(z)} - H_\theta(z) \right] \right)^2 dG(z),$$

$$Q_{n2}(\theta) := \int \left( [\hat{H}_n(z) + H(z)] - \left[ \frac{\mu_n(z, \theta)}{\hat{f}_{Zw}(z)} + H_\theta(z) \right] \right)^2 dG(z).$$

But $Q_{n1}(\theta)$ is bounded above by $2(\rho(\hat{H}_n, H) + \Delta_n(\theta))$, where $\Delta_n(\theta)$ is the $\Delta_n$ of (3.2), with $H$ replaced by $H_\theta$. By (3.5), $\Delta_n(\theta) = o_p(1)$, for each $\theta \in \Theta$.

Similarly, $\forall \theta_1, \theta_2 \in \Theta$, $|\Delta_n(\theta_1) - \Delta_n(\theta_2)|^2$ is bounded above by the product

$$\int \left( \left[ \frac{\mu_n(z, \theta_1)}{\hat{f}_{Zw}(z)} - \frac{\mu_n(z, \theta_2)}{\hat{f}_{Zw}(z)} \right] - [H_{\theta_1}(z) - H_{\theta_2}(z)] \right)^2 dG(z)$$

$$\times \int \left( \left[ \frac{\mu_n(z, \theta_1)}{\hat{f}_{Zw}(z)} + \frac{\mu_n(z, \theta_2)}{\hat{f}_{Zw}(z)} \right] - [H_{\theta_1}(z) + H_{\theta_2}(z)] \right)^2 dG(z).$$

By (m3) and (2.2), the first term of this product is bounded above by $\|\theta_1 - \theta_2\|^{2\beta} O_p(1)$, while the second term is $O_p(1)$ by the boundedness of $m_\theta(x)$ on $\mathcal{I} \times \Theta$. These facts, together with the compactness of $\Theta$, imply that $\sup_{\theta \in \Theta} Q_{n1}(\theta) = o_p(1)$ while $m_\theta(x)$ bounded on $\mathcal{I} \times \Theta$ implies that $\sup_{\theta \in \Theta} Q_{n2}(\theta) = O_p(1)$, thereby completing the proof of (3.6). □

Upon taking $m = m_{\theta_0}$ in the above lemma one immediately obtains the following.

COROLLARY 3.1. *Suppose $H_0$, (e1), (e2), (f1) and (m1)–(m3) hold. Then $\theta_n^* \to \theta_0$, $\hat{\theta}_n \to \theta_0$, in probability.*

**4. Asymptotic distribution of $\hat{\theta}_n$ and $\widehat{\mathcal{D}}_n$.** In this section, we sketch a proof of the asymptotic normality of $\sqrt{n}(\hat{\theta}_n - \theta_0)$ and $\widehat{\mathcal{D}}_n$, under $\mathcal{H}_0$. This proof is similar to that given in [14]. We indicate only the differences. To begin with we focus on $\hat{\theta}_n$. The first step toward this goal is to show that

(4.1) $$nh^d \|\hat{\theta}_n - \theta_0\|^2 = O_p(1).$$

Let $D_n(\theta) = \int \mathcal{Z}_n^2(z, \theta) d\hat{\psi}(z)$. Arguing as in K–N, one obtains

(4.2) $$nh^d D_n(\hat{\theta}_n) = O_p(1).$$

Next, we shall show that for any $a > 0$, there exists an $N_a$ such that

$$(4.3) \quad P\left(D_n(\hat{\theta}_n)/\|\hat{\theta}_n - \theta_0\|^2 \geq a + \inf_{\|b\|=1} b^T \Sigma_0 b\right) > 1 - a \quad \forall n > N_a,$$

where $\Sigma_0$ is as in (m6). The claim (4.1) then follows from (4.2), (4.3), (m6) and the fact $nh^d D_n(\hat{\theta}_n) = nh^d \|\hat{\theta}_n - \theta_0\|^2 [D_n(\hat{\theta}_n)/\|\hat{\theta}_n - \theta_0\|^2]$.

To prove (4.3), let $\Sigma_n(b) := \int [b' \dot{\mu}_n(z, \theta_0)]^2 \, d\hat{\psi}(z)$, $b \in \mathbb{R}^q$ and

$$u_n := \hat{\theta}_n - \theta_0, \qquad d_{ni} := H_{\hat{\theta}_n}(Z_i) - H_{\theta_0}(Z_i) - u_n' \dot{H}_{\theta_0}(Z_i),$$

$$1 \leq i \leq n,$$

(4.4)
$$D_{n1} := \int \left[\frac{1}{n} \sum_{i=1}^n K_{hi}(z) \left(\frac{d_{ni}}{\|u_n\|}\right)\right]^2 d\hat{\psi}(z),$$

$$D_{n2} := \int \left[\frac{u_n'}{\|u_n\|} \dot{\mu}_n(z, \theta_0)\right]^2 d\hat{\psi}(z).$$

Note that

$$\frac{D_n(\hat{\theta}_n)}{\|\hat{\theta}_n - \theta_0\|^2} = \int \frac{\mathcal{Z}_n^2(z, \hat{\theta}_n)}{\|u_n\|^2} d\hat{\psi}(z) \geq D_{n1} + D_{n2} - 2 D_{n1}^{1/2} D_{n2}^{1/2}.$$

We remark here that this inequality corrects a typo in [14] in the equation just above (4.8) on page 120. Assumption (m4) and consistency of $\hat{\theta}_n$ imply that $D_{n1} = o_p(1)$. Exactly the same argument as in [14] with obvious modifications proves that $\sup_{\|b\|=1} \|\Sigma_n(b) - b' \Sigma_0 b\| = o_p(1)$ and (4.3), thereby concluding the proof of (4.1). As in [14], this is used to prove the following theorem where

$$\Sigma = \int \frac{(\sigma_\varepsilon^2 + \tau^2(u)) \dot{H}_{\theta_0}(u) \dot{H}_{\theta_0}'(u) g^2(u)}{f_Z(u)} \, du.$$

THEOREM 4.1. *Assume* (e1)–(e3), (f1), (f2), (g), (k), (m1)–(m5) *and* (h3) *hold. Then under* $\mathcal{H}_0$, $n^{1/2}(\hat{\theta}_n - \theta_0) = \Sigma_0^{-1} n^{1/2} S_n + o_p(1)$. *Consequently,* $n^{1/2} \times (\hat{\theta}_n - \theta_0) \to_d N_q(0, \Sigma_0^{-1} \Sigma \Sigma_0^{-1})$, *where* $\Sigma_0$ *are defined in* (m6).

This theorem shows that asymptotic variance of $n^{1/2}(\hat{\theta}_n - \theta_0)$ consists of two parts. The part involving $\sigma_\varepsilon^2$ reflects the variation in the regression model, while the part involving $\tau^2$ reflects the variation in the measurement error. This is the major difference between asymptotic distribution of the MD estimators discussed for the classical regression model in the K–N paper and for the Berkson model here. Clearly, the larger the measurement error, the larger $\tau^2$ will be.

Next, we state the asymptotic normality result about $\widehat{\mathcal{D}}_n$. Its proof is similar to that of Theorem 5.1 in [14] with obvious modifications and hence no details are given. Recall the notation in (1.4).

THEOREM 4.2. *Suppose* (e1), (e2), (e4), (f1), (f2), (g), (k), (m1)–(m5) *and* (h3) *hold. Then under* $\mathcal{H}_0$, $\widehat{\mathcal{D}}_n \to_d \mathcal{N}_1(0, \Gamma)$ *and* $|\hat{\Gamma}_n \Gamma^{-1} - 1| = o_p(1)$.

Consequently, the test that rejects $\mathcal{H}_0$ whenever $|\widehat{\mathcal{D}}_n| > z_{\alpha/2}$ is of the asymptotic size $\alpha$, where $z_\alpha$ is the $100(1-\alpha)\%$ percentile of the standard normal distribution.

**5. Power of the MD-test.** We shall now discuss some theoretical results about asymptotic power of the proposed tests. We shall show, under some regularity conditions, that $|\widehat{\mathcal{D}}_n| \to \infty$, in probability, under certain fixed alternatives. This in turn implies consistency of the test that rejects $\mathcal{H}_0$ whenever $|\widehat{\mathcal{D}}_n|$ is large against these alternatives. We shall also discuss asymptotic power of the proposed tests against certain local alternatives. Accordingly, let $m \in L_2(G)$ and $H(z) := E(m(X)|Z = z)$. Also, let $\nu(z, \theta) := H_\theta(z) - H(z)$, $e_{ni} := Y_i - H_{\theta_n}(Z_i)$, $e_i := Y_i - H(Z_i)$, where $\theta_n$ is an estimator of $T(H)$ of (1.5). Let, for $z \in \mathbb{R}^d, \theta \in \Theta$,

$$\mathcal{V}_n(z) := \frac{1}{n}\sum_{i=1}^n K_{hi}(z)e_i, \qquad \bar{\nu}_n(z, \theta) := \frac{1}{n}\sum_{i=1}^n K_{hi}(z)\nu(Z_i, \theta).$$

5.1. *Consistency.* Let $\mathcal{D}_n := nh^{d/2}\mathcal{G}_n^{-1/2}(M_n(\theta_n) - C_n)$, where

$$C_n := \frac{1}{n^2}\sum_{i=1}^n \int K_{hi}^2 e_{ni}^2 \, d\hat{\psi}, \qquad \mathcal{G}_n := 2n^{-2}h^d \sum_{i \neq j}\left(\int K_{hi}K_{hj}e_{ni}e_{nj} \, d\hat{\psi}\right)^2.$$

If $\theta_n = \hat{\theta}_n$, then $C_n = \hat{C}_n, \mathcal{G}_n = \hat{\Gamma}_n$ and $\mathcal{D}_n = \widehat{\mathcal{D}}_n$. The following theorem provides a set of sufficient conditions under which $|\mathcal{D}_n| \to \infty$, in probability, for any sequence of consistent estimator $\theta_n$ of $T(H)$.

THEOREM 5.1. *Suppose* (e1), (e2), (e4), (f1), (f2), (g), (k), (m3), (h3) *and the alternative hypothesis* $\mathcal{H}_1 : \mu(x) = m(x), \forall x \in \mathcal{I}$ *hold with the additional assumption that* $\inf_\theta \rho(H, H_\theta) > 0$. *Then, for any sequence of consistent estimator* $\theta_n$ *of* $T(H)$, $|\mathcal{D}_n| \to \infty$, *in probability. Consequently,* $|\widehat{\mathcal{D}}_n| \to \infty$, *in probability.*

PROOF. Subtracting and adding $H(Z_i)$ from $e_{ni}$, we obtain $M_n(\theta_n) = S_{n1} - 2S_{n2} + S_{n3}$, where $S_{n1} := \int \mathcal{V}_n^2 \, d\hat{\psi}$, $S_{n2} := \int \mathcal{V}_n(z)\bar{\nu}_n(z, \theta_n) \, d\hat{\psi}(z)$ and $S_{n3} := \int \bar{\nu}_n^2(z, \theta_n) \, d\hat{\psi}(z)$. Arguing as in Lemma 5.1 of [14], we can verify that under the current setup, $nh^{d/2}(S_{n1} - C_n^*) \to_d N_1(0, \Gamma^*)$, where $C_n^* = \sum_{i=1}^n \int K_{hi}^2(z) e_i^2 \, d\hat{\psi}(z)/n^2$, $\Gamma^* = 2 \int (\sigma_*^2(z))^2 g(z) \, d\psi(z) \|K_2\|^2$ with $\sigma_*^2(z) = \sigma_e^2 + \gamma^2(z)$, with $\sigma_e^2(z) = E[(Y - H(Z))^2 | Z = z]$.

Next, consider $S_{n3}$. For convenience write $T$ for $T(H)$. By subtracting and adding $H_T(Z_i)$ from $\nu(Z_i, \theta_n)$, we have $S_{n3} = S_{n31} + 2S_{n32} + S_{n33}$, where

$$S_{n31} := \int \bar{\nu}_n^2(z, T) \, d\hat{\psi}(z),$$

$$S_{n32} := \int \bar{\nu}_n(z, T)[\bar{\mu}_n(z, \theta_n) - \bar{\mu}_n(z, T)] \, d\hat{\psi}(z),$$

$$S_{n33} := \int [\bar{\mu}_n(z, \theta_n) - \bar{\mu}_n(z, T)]^2 \, d\hat{\psi}(z).$$

Routine calculations and (2.2) show that $S_{n31} = \rho(H, H_T) + o_p(1)$, under $\mathcal{H}_1$. By (m3), $S_{n33} \leq \|\theta_n - T\|^{2\beta} \int_\mathcal{I} [\frac{1}{n \hat{f}_{Zw}(z)} \sum_{i=1}^n K_{hi}(z) |l(Z_i)|]^2 \, dG(z) = o_p(1)$, by consistency of $\theta_n$ for $T$. By the C–S inequality, one obtains that $S_{n32} = o_p(1) = S_{n2}$. Therefore, $S_{n3} = \rho(H, H_T) + o_p(1)$.

Note that

$$C_n - C_n^* = -\frac{2}{n^2} \sum_{i=1}^n \int K_{hi}^2(z) e_i \nu(Z_i, \theta_n) \, d\hat{\psi}(z)$$

$$+ \frac{1}{n^2} \sum_{i=1}^n \int K_{hi}^2(z) \nu^2(Z_i, \theta_n) \, d\hat{\psi}(z).$$

Both terms on the right-hand side are of the order $o_p(1)$.

We shall next show that $\mathcal{G}_n \to \Gamma^*$ in probability. Adding and subtracting $H(Z_i)$ and $H(Z_j)$ from $e_{ni}$ and $e_{nj}$, respectively, and expanding the square of integral, one can rewrite $\mathcal{G}_n = \sum_{j=1}^{10} A_{nj}$, where

$$A_{n1} = 2h^d n^{-2} \sum_{i \neq j} \left( \int K_{hi}(z) K_{hj}(z) e_i e_j \, d\hat{\psi}(z) \right)^2,$$

$$A_{n2} = 2h^d n^{-2} \sum_{i \neq j} \left( \int K_{hi}(z) K_{hj}(z) e_i \nu(Z_j, \theta_n) \, d\hat{\psi}(z) \right)^2,$$

$$A_{n3} = 2h^d n^{-2} \sum_{i \neq j} \left( \int K_{hi}(z) K_{hj}(z) \nu(Z_i, \theta_n) e_j \, d\hat{\psi}(z) \right)^2,$$

$$A_{n4} = 2h^d n^{-2} \sum_{i \neq j} \left( \int K_{hi}(z) K_{hj}(z) \nu(Z_i, \theta_n) \nu(Z_j, \theta_n) \, d\hat{\psi}(z) \right)^2,$$

$$A_{n5} = -4h^d n^{-2} \sum_{i \neq j} \left( \int K_{hi}(z) K_{hj}(z) e_i e_j \, d\hat{\psi}(z) \right.$$

$$\times \int K_{hi}(z) K_{hj}(z) e_i \nu(Z_j, \theta_n) \, d\hat{\psi}(z) \Big),$$

$$A_{n6} = -4h^d n^{-2} \sum_{i \neq j} \Big( \int K_{hi}(z) K_{hj}(z) e_i e_j \, d\hat{\psi}(z)$$

$$\times \int K_{hi}(z) K_{hj}(z) \nu(Z_i, \theta_n) e_j \, d\hat{\psi}(z) \Big),$$

$$A_{n7} = 4h^d n^{-2} \sum_{i \neq j} \Big( \int K_{hi}(z) K_{hj}(z) e_i e_j \, d\hat{\psi}(z)$$

$$\times \int K_{hi}(z) K_{hj}(z) \nu(Z_i, \theta_n) \nu(Z_j, \theta_n) \, d\hat{\psi}(z) \Big),$$

$$A_{n8} = 4h^d n^{-2} \sum_{i \neq j} \Big( \int K_{hi}(z) K_{hj}(z) e_i \nu(Z_j, \theta_n) \, d\hat{\psi}(z)$$

$$\times \int K_{hi}(z) K_{hj}(z) \nu(Z_i, \theta_n) e_j \, d\hat{\psi}(z) \Big),$$

$$A_{n9} = -4h^d n^{-2} \sum_{i \neq j} \Big( \int K_{hi}(z) K_{hj}(z) e_i \nu(Z_j, \theta_n) \, d\hat{\psi}(z)$$

$$\times \int K_{hi}(z) K_{hj}(z) \nu(Z_i, , \theta_n) \nu(Z_j, \theta_n) \, d\hat{\psi}(z) \Big),$$

$$A_{n10} = -4h^d n^{-2} \sum_{i \neq j} \Big( \int K_{hi}(z) K_{hj}(z) \nu(Z_i, \theta_n) e_j \, d\hat{\psi}(z)$$

$$\times \int K_{hi}(z) K_{hj}(z) \nu(Z_i, \theta_n) \nu(Z_j, \theta_n) \, d\hat{\psi}(z) \Big).$$

By taking the expectation, using Fubini's theorem we obtain

(5.1) $\quad h^d n^{-2} \sum_{i \neq j} \Big( \int K_{hi}(z) K_{hj}(z) |e_i||e_j| \, d\psi(z) \Big)^2 = O_p(1),$

(5.2) $\quad h^d n^{-2} \sum_{i \neq j} \Big( \int K_{hi}(z) K_{hj}(z) |e_i|^k \, d\psi(z) \Big)^2 = O_p(1), \qquad k = 0, 1.$

By (2.2) and (5.1) and arguing as in the proof of Lemma 5.5 in K–N, one can verify that $A_{n1} \to_p \Gamma_1^* := 2 \int (\sigma_e^4 \frac{g^2}{f_Z^2})(z) \, dz \|K_2\|^2$.

Add and subtract $H_T(Z_j)$ from $\nu(Z_j, \theta_n)$, to obtain

$$A_{n2} = \frac{2h^d}{n^2} \sum_{i \neq j} \Big( \int K_{hi}(z) K_{hj}(z) e_i \nu(Z_j, \theta) \, d\hat{\psi}(z) \Big)^2$$

$$+ \frac{4h^d}{n^2}\sum_{i\neq j}\biggl(\int K_{hi}(z)K_{hj}(z)e_i\nu(Z_j,\theta)\,d\hat{\psi}(z)$$

$$\times \int K_{hi}(z)K_{hj}(z)e_i(H_{\theta_n}(Z_j) - H_T(Z_j))\,d\hat{\psi}(z)\biggr)$$

$$+ \frac{2h^d}{n^2}\sum_{i\neq j}\biggl(\int K_{hi}(z)K_{hj}(z)e_i(H_{\theta_n}(Z_j) - H_T(Z_j))\,d\hat{\psi}(z)\biggr)^2.$$

By (m4), consistency of $\theta_n$, (2.1), the C–S inequality on the double sum and (5.2), the last two terms of the above expression are $o_p(1)$. Arguing, as for $A_{n1}$, the first term on the right-hand side above converges in probability to $\Gamma_2^* := 2\int\sigma_e^2(z)[H(z) - H_T(z)]^2\frac{g^2(z)}{f_Z^2(z)}\,dz\|K_2\|^2$. Similarly, one can also show $A_{n3} \to \Gamma_2^*$ in probability.

Similarly, by adding and subtracting $H_T(Z_i)$, $H_T(Z_j)$ from $\nu(Z_i,\theta_n)$, $\nu(Z_j,\theta_n)$, respectively, in $A_{n4}$, one obtains $A_{n4} = \Gamma_3^* + o_p(1)$, where $\Gamma_3^* = 2\int[H(z) - H_T(z)]^4\frac{g^2(z)}{f_Z^2(z)}\,dz\|K_2\|^2$. Next, rewrite

$$A_{n5} = -4h^d n^{-2}\sum_{i\neq j}\biggl(\int K_{hi}(z)K_{hj}(z)e_ie_j\,d\hat{\psi}(z)$$

$$\times \int K_{hi}(z)K_{hj}(z)e_i\nu(Z_j,\theta_n)\,d\hat{\psi}(z)\biggr)$$

$$- 4h^d n^{-2}\sum_{i\neq j}\biggl(\int K_{hi}(z)K_{hj}(z)e_ie_j\,d\hat{\psi}(z)$$

$$\times \int K_{hi}(z)K_{hj}(z)e_i[H_{\theta_n}(Z_j) - H_T(Z_j)]\,d\hat{\psi}(z)\biggr)$$

$$= A_{n51} + A_{n52}, \quad\text{say.}$$

Clearly, $E\tilde{A}_{n51} = 0$. Argue as for (5.13) in K–N, verify that $E(\tilde{A}_{n51}^2) = O((nd)^{-1})$. Therefore, $\tilde{A}_{n51} = o_p(1)$. By the C–S inequality on the double sum, (2.2), (5.1) and (5.2), we have $A_{n51} = \tilde{A}_{n51} + o_p(1)$. Hence $A_{n51} = o_p(1)$. Similarly, one can verify $A_{n52} = o_p(1)$. These results imply $A_{n5} = o_p(1)$.

Similarly, one can show that $A_{ni} = o_p(1)$, $i = 6, 7, 8, 9, 10$. Note that $\Gamma^* = \Gamma_1^* + 2\Gamma_2^* + \Gamma_3^*$, so we obtain that $\mathcal{G}_n \to \Gamma^*$, in probability.

All these results together imply that

$$\mathcal{D}_n = nh^{d/2}\hat{\Gamma}_n^{-1/2}(S_{n1} - C_n^*) + nh^{d/2}\mathcal{G}_n^{-1/2}\rho(H, H_T) + o_p(nh^{d/2}),$$

hence the theorem. $\square$

5.2. *Power at local alternatives.* Here we shall now study the asymptotic power of the proposed MD-test against some local alternatives. Accordingly, let $r$ be a known continuously differentiable real-valued function and let $R(z) := E(r(X)|Z=z)$. In addition, assume $R \in L_2(G)$ and

$$\int H_\theta R \, dG = 0 \qquad \forall \theta \in \Theta. \tag{5.3}$$

Consider the sequence of local alternatives

$$\mathcal{H}_{1n} : \mu(x) = m_{\theta_0}(x) + \gamma_n r(x), \qquad \gamma_n = 1/\sqrt{nh^{d/2}}. \tag{5.4}$$

The following theorem gives the asymptotic distribution of $\hat{\theta}_n$ under $\mathcal{H}_{1n}$.

THEOREM 5.2. *Suppose* (e1)–(e3), (f1), (f2), (g), (k), (m1)–(m5) *and* (h3) *hold; then under the local alternative (5.3) and (5.4)*, $n^{1/2}(\hat{\theta}_n - \theta_0) \to_d N_q(0, \Sigma_0^{-1} \Sigma \Sigma_0^{-1})$.

PROOF. The basic idea of the proof is the same as in the null case. We only stress the differences here. Under $\mathcal{H}_{1n}$, $\varepsilon_i \equiv Y_i - m_{\theta_0}(X_i) - \gamma_n r(X_i)$. Let $\bar{r}_n(z) := \sum_{i=1}^n K_{hi}(z) r(X_i)/n$.

We first note that $nh^d M_n(\theta_0) = O_p(1)$. In fact, under (5.4), $M_n(\theta_0)$ can be bounded above by 2 times the sum of $(1/nh^{d/2}) \int \bar{r}_n^2 \, d\hat{\psi}$ and

$$\int \left[ \sum_{i=1}^n K_{hi}(z)(m_{\theta_0}(X_i) + \varepsilon_i - H_{\theta_0}(Z_i)) \right]^2 d\hat{\psi}(z).$$

Using the variance argument and (2.2), one verifies that this term is of the order $O_p(n^{-1}h^{-d})$. Note that $\bar{r}_n$ is a kernel estimator of $R$. Hence, $R \in L_2(G)$ and a routine argument shows that the former term is $O_p(n^{-1}h^{-d/2})$. This leads to the conclusion $nh^d M_n(\theta_0) = O_p(1)$. This fact and an argument similar to the one used in K–N, together with the fact $\hat{\theta}_n \to_p \theta_0$, yield $nh^d \|\hat{\theta}_n - \theta_0\| = O_p(1)$, under $\mathcal{H}_{1n}$.

Note that with $\dot{M}_n(\theta) := \partial M_n(\theta)/\partial \theta$, $\hat{\theta}_n$ satisfies

$$\dot{M}_n(\hat{\theta}_n) = -2 \int U_n(z, \hat{\theta}_n) \dot{\mu}_n(z, \hat{\theta}_n) \, d\hat{\psi}(z) = 0, \tag{5.5}$$

where $U_n(z, \theta)$ and $\dot{\mu}_n(z, \theta)$ are defined in (2.4). Adding and subtracting $H_{\theta_0}(Z_i)$ from $Y_i - H_{\hat{\theta}_n}(Z_i)$ in $U_n(z, \hat{\theta}_n)$, we can rewrite (5.5) as

$$\int U_n(z) \dot{\mu}_n(z, \hat{\theta}_n) \, d\hat{\psi}(z) = \int \mathcal{Z}_n(z, \hat{\theta}_n) \dot{\mu}_n(z, \hat{\theta}_n) \, d\hat{\psi}(z). \tag{5.6}$$

The right-hand side of (5.6) involves the error variables only through $\hat{\theta}_n$. Since under $\mathcal{H}_{1n}$ we also have $n^{1/2}(\hat{\theta}_n - \theta_0) = O_p(1)$, its asymptotic behavior under $\mathcal{H}_{1n}$ is the same as in the null case, that is, it equals $\mathcal{R}_n(\hat{\theta}_n - \theta_0) +$ $o_P(1), \mathcal{R}_n = \Sigma_0 + o_p(1)$. The left-hand side, under (5.4), can be rewritten as $\mathcal{S}_{n1} + \mathcal{S}_{n2}$, where

$$\mathcal{S}_{n1} = \int \frac{1}{n} \sum_{i=1}^{n} K_{hi}(z)[m_{\theta_0}(X_i) + \varepsilon_i - H_{\theta_0}(Z_i)] \dot{\mu}_n(z, \hat{\theta}_n) \, d\hat{\psi}(z),$$

$$\mathcal{S}_{n2} = \gamma_n \int \bar{r}_n(z) \dot{\mu}_n(z, \hat{\theta}_n) \, d\hat{\psi}(z).$$

Note that $m_{\theta_0}(X_i) + \varepsilon_i - H_{\theta_0}(Z_i)$ are i.i.d. with mean 0 and finite second moment. Arguing as in the proofs of Lemmas 4.1 and 4.2 of [14] with $\varepsilon_i$ there replaced by $m_{\theta_0}(X_i) + \varepsilon_i - H_{\theta_0}(Z_i)$ yields that under $\mathcal{H}_{1n}$, $\sqrt{n}\mathcal{S}_{n1} \to_d N_q(0, \Sigma)$. Thus, the theorem will be proved if we can show $\sqrt{n}\mathcal{S}_{n2} = o_p(1)$.

For this purpose, with $r_h(z) := EK_h(z - Z)r(X) = EK_h(z - Z)R(Z)$, we need the following facts. Arguing as for (3.4) and using differentiability of $r$, one obtains

(5.7) $$\int [\bar{r}_n - r_h]^2 \, d\hat{\psi} = O_p(n^{-1}h^{-d}), \qquad \int [r_h - Rf_Z]^2 \, d\psi = O(h^{2d}),$$

$$\int \|\dot{\mu}_h - \dot{H}_{\theta_0} f_Z\|^2 \, d\psi = O(h^{2d}).$$

Then the integral in $\mathcal{S}_{n2}$ can be written as

$$\int \{[\bar{r}_n(z) - r_h(z)] + [r_h(z) - R(z)f_Z(z)] + R(z)f_Z(z)\}$$
$$\times \{[\dot{\mu}_n(z, \hat{\theta}_n) - \dot{\mu}_n(z, \theta_0)] + [\dot{\mu}_n(z, \theta_0) - \dot{\mu}_h(z)]$$
$$+ [\dot{\mu}_h(z) - \dot{H}_{\theta_0}(z)f_Z(z)] + \dot{H}_{\theta_0}(z)f_Z(z)\} \, d\hat{\psi}(z).$$

This can be further expanded into twelve terms. By (m5), (2.2), (5.7) and C–S, one can show that all of these twelve terms are $o_p(h^{-d/4})$ except the term $\int R\dot{H}_{\theta_0} f_Z^2 \, d\hat{\psi} = \int R\dot{H}_{\theta_0} \, dG + \int R\dot{H}_{\theta_0} q_n \, dG$. But (5.3), (m6), continuity of $r(x)$ and the compactness of $\mathcal{I}$ imply $\int \dot{H}_{\theta_0} R \, dG = 0$. The second term is bounded above by

(5.8) $$\sup_{z \in \mathcal{I}} |q_n(z)| \int |R| \|\dot{H}_{\theta_0}\| \, dG.$$

By Theorem 2.2, part (2) in Bosq [4] and the choice of $w = (\frac{\log n}{n})^{1/(d+4)}$, $(\log_k n)^{-1}(n/\log n)^{2/(d+4)} \sup_{z \in \mathcal{I}} |\hat{f}_{Zw}(z) - f_Z(z)| \to 0$, almost surely, for all $k > 0$. This fact and (h3) readily imply that (5.8) is of the order $o_p(h^{d/4})$, so that $n^{1/2}\mathcal{S}_{n2} = \sqrt{n} \cdot (\sqrt{n}h^{d/2})^{-1} \cdot o_p(h^{d/4}) = o_p(1)$. Hence the theorem. $\square$

The following theorem gives asymptotic power of the MD-test against the local alternative (5.3) and (5.4).

THEOREM 5.3. *Suppose* (e1), (e2), (e4), (f1), (f2), (g), (k), (m4), (h3) *and the local alternative hypothesis (5.3) and (5.4) hold. Then,* $\widehat{\mathcal{D}}_n \to_d N(\Gamma^{-1/2} \int R^2 \, dG, 1)$, *where* $\Gamma$ *is as in (2.3).*

PROOF. Rewrite $M_n(\hat{\theta}_n) = T_{n1} + 2T_{n2} + T_{n3}$, where $T_{n1} := \int U_n^2 \, d\hat{\psi}$, $T_{n2} := \int U_n(z)[\mu_n(z,\theta_0) - \mu_n(z,\hat{\theta}_n)] \, d\hat{\psi}$ and $T_{n3} := \int [\mu_n(z,\theta_0) - \mu_n(z,\hat{\theta}_n)]^2 \, d\hat{\psi}$. By Theorem 5.2, $\sqrt{n}(\hat{\theta}_n - \theta_0) = O_p(1)$. This fact, (m4) and (2.2) imply $T_{n3} = O_p(n^{-1})$.

Next, we shall show that $T_{n2} = O_p(n^{-1}h^{-d/4})$. By C–S, $T_{n2}^2 \leq T_{n1}T_{n3}$. Moreover, $T_{n1} = \int U_n^2 \, d\psi + \int U_n^2 q_n \, d\psi$. But under $\mathcal{H}_{1n}$, $Y_i = m_{\theta_0}(X_i) + \gamma_n r(X_i) + \varepsilon_i$. Hence, $\int U_n^2 \, d\psi$ is bounded above by 3 times the sum

$$\int \left[\frac{1}{nf_Z(z)} \sum_{i=1}^n K_{hi}(z)\varepsilon_i\right]^2 dG(z) + \int \bar{r}_n^2 \, d\hat{\psi},$$

$$+ \int \left[\frac{1}{nf_Z(z)} \sum_{i=1}^n K_{hi}(z)[m_{\theta_0}(X_i) - H_{\theta_0}(Z_i)]\right]^2 dG(z).$$

Arguing as in Section 2, all of these terms are $O_p(n^{-1}h^{-d/2})$. This fact and (2.2) imply that the second term in $T_{n1}$ is of the order $o_p(n^{-1}h^{-d/2})$. Hence $T_{n1} = O_p(n^{-1}h^{-d/2})$ and $T_{n2} = O_p(n^{-1}h^{-d/4})$.

We shall now obtain a more precise approximation to $T_{n1}$. For this purpose, write $\xi_i = \varepsilon_i + m_{\theta_0}(X_i) - H_{\theta_0}(Z_i)$ and let $V_n(z) := \sum_{i=1}^n K_{hi}(z)\xi_i/n$. Then, $T_{n1} = T_{n11} + 2\gamma_n T_{n12} + \gamma_n^2 T_{n13}$, where

$$T_{n11} := \int V_n \, d\hat{\psi}, \qquad T_{n12} := \int V_n \bar{r}_n \, d\hat{\psi}, \qquad T_{n13} := \int \bar{r}_n^2 \, d\hat{\psi}.$$

Now, we shall show that

$$(5.9) \qquad \int (R/\hat{f}_{Zw}) V_n \, dG = o_p(1/\sqrt{nh^{d/2}}).$$

In fact, with $d\psi_1 := dG/f_Z$, the left-hand side equals $\int (R/f_Z) V_n \, dG + \int R V_n q_{w1} \, d\psi_1$. The first term is an average of i.i.d. mean-zero r.v.'s and a variance calculation shows that it is of the order $O_p(n^{-1/2})$, while by Theorem 2.2, part (2) in Bosq [4], the second term is of the order $o_p(1/\sqrt{nh^{d/2}})$, thereby proving (5.9).

Arguing as for (3.4) one obtains that $\int V_n^2 \, d\psi = O_p(1/nh^d)$. Next, note that $\bar{r}_n/\hat{f}_{Zw}$ is an estimator of $R$, so by the C–S inequality again,

$$\int (V_n/\hat{f}_{Zw})[(\bar{r}_n/\hat{f}_{Zw}) - R] \, dG = o_p(1/\sqrt{nh^{d/2}}).$$

This fact and (5.9) imply that $T_{n12} = o_p(1/\sqrt{nh^{d/2}})$. A similar and relatively easier argument yields that $T_{n13} = \int R^2 \, dG + o_p(1)$.

Finally, we need to discuss asymptotic behavior of $C_n$ under the local alternative (5.4). With $\zeta_i = Y_i - H_{\theta_0}(Z_i)$, rewrite $Y_i - H_{\theta_n}(Z_i) = \zeta_i + H_{\theta_0}(Z_i) - H_{\theta_n}(Z_i)$ in $C_n$, to obtain

$$C_n = \frac{1}{n^2} \sum_{i=1}^n \int K_{hi}^2(z) \zeta_i^2 \, d\hat{\psi}(z)$$
$$+ \frac{2}{n^2} \sum_{i=1}^n \int K_{hi}^2(z) \zeta_i (H_{\theta_0}(Z_i) - H_{\theta_n}(Z_i)) \, d\hat{\psi}(z)$$
$$+ \frac{1}{n^2} \sum_{i=1}^n \int K_{hi}^2(z) (H_{\theta_0}(Z_i) - H_{\theta_n}(Z_i))^2 \, d\hat{\psi}(z)$$
$$= C_{n1} + 2C_{n2} + C_{n3}.$$

But with notation at (4.4),

$$C_{n2} = -\frac{1}{n^2} \sum_{i=1}^n \int K_{hi}^2(z) \xi_i d_{ni} \, d\hat{\psi}(z) - \frac{1}{n^2} \sum_{i=1}^n \int K_{hi}^2(z) r(X_i) d_{ni} \, d\hat{\psi}(z)$$
$$+ \frac{1}{n^2} \sum_{i=1}^n \int K_{hi}^2(z) \xi_i u_n' \dot{H}_{\theta_0}(Z_i) \, d\hat{\psi}(z)$$
$$+ \frac{1}{n^2} \sum_{i=1}^n \int K_{hi}^2(z) r(X_i) u_n' \dot{H}_{\theta_0}(Z_i) \, d\hat{\psi}(z).$$

Recall that $\gamma_n = 1/\sqrt{nh^{d/2}}$. Using assumptions (m4), (h2), one can show the first and the third terms in $C_{n2}$ are of the order $O_P(n^{-3/2}h^{-d})$, the second and the fourth terms are of the order $O_p(n^{-2}h^{-3d/2})$. This implies $C_{n2} = o_p(\gamma_n^2)$. Similarly, one can show that $C_{n3} = O_p(n^{-3/2}h^{-d}) = o_p(\gamma_n^2)$.

Since $Y_i - H_{\theta_0}(Z_i) = \xi_i + \gamma_n r(X_i)$, if we let $D_n = n^{-2} \sum_{i=1}^n \int K_{hi}^2 \xi_i^2 \, d\hat{\psi}$, then using the similar argument, we can show that $C_{n1} = D_n + o_p(\gamma_n^2)$.

To see the asymptotic property of $\hat{\Gamma}_n$ under the local alternative, adding and subtracting $H_{\theta_0}(Z_i)$, $H_{\theta_0}(Z_j)$ from $e_{ni}$ and $e_{ni}$, respectively, and letting $\xi_i = m_{\theta_0}(X_i) - H_{\theta_0}(Z_i) + \varepsilon_i$, we will arrive at

$$\hat{\Gamma}_n = 2n^{-2}h^d \sum_{i \neq j} \left( \int K_{hi}(z) K_{hj}(z) \xi_i \xi_j \, d\hat{\psi}(z) \right)^2 + \omega_n.$$

The first term converges in probability to $\Gamma$. The remainder $\omega_n = o_p(1)$ can be proven by using the C–S inequality on the double sum, consistency of $\theta_n$, (2.2) and the following facts:

$$\frac{h^d}{n^2} \sum_{i \neq j} \left( \int K_{hi}(z) K_{hj}(z) |\xi_i| |\xi_j| \, d\psi(z) \right)^2 = O_p(1),$$

$$\frac{h^d}{n^2}\sum_{i\neq j}\left(\int K_{hi}(z)K_{hj}(z)|\xi_i|^k\,d\psi(z)\right)^2 = O_p(1), \qquad k=0,1.$$

Therefore, under the local alternative hypothesis (5.4),

$$nh^{d/2}\hat{\Gamma}_n^{-1/2}(M_n(\hat{\theta}_n) - \hat{C}_n) = nh^{d/2}\hat{\Gamma}_n^{-1/2}(T_{n11} - D_n) + \hat{\Gamma}_n^{-1/2}T_{n13} + o_p(1),$$

which, together with the fact $nh^{d/2}(T_{n11} - D_n) \to_d N_1(0,\Gamma)$, $T_{n13} \to \int R^2\,dG$ and $\hat{\Gamma}_n \to \Gamma$ in probability, implies the theorem. $\square$

**6. Simulations.** This section contains results of two simulation studies corresponding to the following cases: case 1: $d=q=1$ and $m_\theta$ linear; case 2: $d=q=2$ and $m_\theta$ nonlinear. In each case the Monte Carlo average values of $\hat{\theta}_n$, $\mathrm{MSE}(\hat{\theta}_n)$, empirical levels and powers of the MD test are reported. The asymptotic level is taken to be 0.05 in all cases.

In the first case $\{Z_i\}_{i=1}^n$ are obtained as a random sample from the uniform distribution on $[-1,1]$ and $\{\varepsilon_i\}_{i=1}^n$ and $\{\eta_i\}_{i=1}^n$ are obtained as two independent random samples from $\mathcal{N}_1(0,(0.1)^2)$. Then $(X_i,Y_i)$ are generated using the model $Y_i = \mu(X_i) + \varepsilon_i$, $X_i = Z_i + \eta_i$, $i=1,2,\ldots,n$.

The kernel functions and the bandwidths used in the simulation are

$$K(z) = K^*(z) = \tfrac{3}{4}(1-z^2)I(|z|\leq 1), \qquad h = \frac{a}{n^{1/3}}, \qquad w = b\left(\frac{\log n}{n}\right)^{1/5},$$

with some choices for $a$ and $b$. The integrating measure $G$ is taken to be the uniform measure on $[-1,1]$.

The parametric model is taken to be $m_\theta(x) = \theta x$, $x,\theta \in \mathbb{R}$, $\theta_0 = 1$. Then, $H_\theta(z) = \theta z$. In this case various calculations simplify as follows. By taking the derivative of $M_n(\theta)$ in $\theta$ and solving the equation of $\partial M_n(\theta)/\partial\theta = 0$, we obtain $\hat{\theta}_n = A_n/B_n$, where

$$A_n = \int_{-1}^{1}\left[\sum_{i=1}^n K_{hi}(z)Y_i\right]\left[\sum_{i=1}^n K_{hi}(z)Z_i\right]\left[\sum_{i=1}^n K_{wi}(z)\right]^{-2}dz,$$

$$B_n = \int_{-1}^{1}\left[\sum_{i=1}^n K_{hi}(z)Z_i\right]^2\left[\sum_{i=1}^n K_{wi}(z)\right]^{-2}dz.$$

Then, with $\hat{\varepsilon}_i := Y_i - \hat{\theta}_n Z_i$,

$$M_n(\hat{\theta}_n) = \int_{-1}^{1}\left(\sum_{i=1}^n K_{hi}(z)\hat{\varepsilon}_i\right)^2\left(\sum_{i=1}^n K_{wi}(z)\right)^{-2}dz,$$

$$\hat{C}_n = \int_{-1}^{1}\left(\sum_{i=1}^n K_{hi}^2(z)\hat{\varepsilon}_i^2\right)\left(\sum_{i=1}^n K_{wi}(z)\right)^{-2}dz.$$

Table 1 reports the Monte Carlo mean and $\mathrm{MSE}(\hat{\theta}_n)$ under $H_0$ for the sample sizes $50, 100, 200, 500$, each repeated 1000 times. One can see there appears to be little bias in $\hat{\theta}_n$ for all chosen sample sizes and as expected, the MSE decreases as the sample size increases.

To assess the level and power behavior of the $\widehat{\mathcal{D}}_n$-test, we chose the following four models to simulate data from; in each of these cases $X_i = Z_i + \eta_i$:

Model 0: $Y_i = X_i + \varepsilon_i$,

Model 1: $Y_i = X_i + 0.3X_i^2 + \varepsilon_i$,

Model 2: $Y_i = X_i + 1.4\exp(-0.2X_i^2) + \varepsilon_i$,

Model 3: $Y_i = X_i I(X_i \geq 0.2) + \varepsilon_i$.

To assess the effect of the choice of $(a,b)$ that appear in the bandwidths on the level and power, we ran simulations for numerous choices of $(a,b)$, ranging from 0.2 to 1. Table 2 reports these simulation results pertaining to $\widehat{\mathcal{D}}_n$ for three choices of $(a,b)$. Simulation results for the other choices were similar to those reported here. Data from Model 0 in this table are used to study empirical sizes and data from Models 1 to 3 are used to study empirical powers of the test. These entities are obtained by computing $\#\{|\widehat{\mathcal{D}}_n| \geq 1.96\}/1000$, where $\widehat{\mathcal{D}}_n := nh^{d/2}\hat{\Gamma}_n^{-1/2}(M_n(\hat{\theta}_n) - \hat{C}_n)$.

From Table 2, one sees that empirical level is sensitive to the choice of $(a,b)$ for moderate sample sizes ($n \leq 200$) but gets closer to the asymptotic level of 0.05 with the increase in the sample size and hence is stable over the chosen values of $(a,b)$ for large sample sizes. On the other hand the empirical power appears to be far less sensitive to the values of $(a,b)$ for the sample sizes of 100 and more. Even though the theory of the present paper is not applicable to Model 3, it was included here to see the effect of the discontinuity in the regression function on the power of the minimum distance test. In our simulation, the discontinuity of the regression has little effect on the power of the minimum distance test.

Now consider the case 2 where $d = 2, q = 2$ and $m_\theta(x) = \theta_1 x_1 + \exp(\theta_2 x_2)$, $\theta = (\theta_1, \theta_2)' \in \mathbb{R}^2$, $x_1, x_2 \in \mathbb{R}$. Accordingly, here $H_\theta(z) = \theta_1 z_1 + \exp(\theta_2 z_2 + 0.005\theta_2^2)$. The true $\theta_0 = (1, 2)'$ was used in these simulations.

In all models below, $\{Z_i = (Z_{1i}, Z_{2i})'\}_{i=1}^n$ are obtained as a random sample from the uniform distribution on $[-1, 1]^2$, $\{\varepsilon_i\}_{i=1}^n$ are obtained from

TABLE 1
*Mean and MSE of $\hat{\theta}_n$*

| Sample size | 50 | 100 | 200 | 500 |
|---|---|---|---|---|
| Mean | 1.0003 | 0.9987 | 1.0006 | 0.9998 |
| MSE | 0.0012 | 0.0006 | 0.0003 | 0.0001 |

TABLE 2
*Levels and powers of the minimum distance test*

| Model | a, b | Sample size | | | |
|---|---|---|---|---|---|
| | | 50 | 100 | 200 | 500 |
| Model 0 | 0.3, 0.2 | 0.007 | 0.026 | 0.028 | 0.048 |
| | 0.5, 0.5 | 0.014 | 0.022 | 0.040 | 0.051 |
| | 1.0, 1.0 | 0.021 | 0.020 | 0.031 | 0.043 |
| Model 1 | 0.3, 0.2 | 0.754 | 0.987 | 1.000 | 1.000 |
| | 0.5, 0.5 | 0.945 | 1.000 | 1.000 | 1.000 |
| | 1.0, 1.0 | 1.000 | 1.000 | 1.000 | 1.000 |
| Model 2 | 0.3, 0.2 | 0.857 | 0.996 | 1.000 | 1.000 |
| | 0.5, 0.5 | 0.999 | 1.000 | 1.000 | 1.000 |
| | 1.0, 1.0 | 1.000 | 1.000 | 1.000 | 1.000 |
| Model 3 | 0.3, 0.2 | 0.874 | 0.993 | 1.000 | 1.000 |
| | 0.5, 0.5 | 1.000 | 1.000 | 1.000 | 1.000 |
| | 1.0, 1.0 | 1.000 | 1.000 | 1.000 | 1.000 |

$\mathcal{N}_1(0,(0.1)^2)$ and $\{\eta_i = (\eta_{1i}, \eta_{2i})'\}_{i=1}^n$ are obtained from the bivariate normal distribution with mean vector 0 and the diagonal covariance matrix with both diagonal entries equal to $(0.1)^2$. We simulated data from the following four models, where $X_i = Z_i + \eta_i$:

Model 0: $Y_i = X_{1i} + \exp(2X_{2i}) + \varepsilon_i$,

Model 1: $Y_i = X_{1i} + \exp(2X_{2i}) + 1.4X_{1i}^2 + 1 + \varepsilon_i$,

Model 2: $Y_i = X_{1i} + \exp(2X_{2i}) + 1.4X_{1i}^2 X_{2i}^2 + \varepsilon_i$,

Model 3: $Y_i = X_{1i} + \exp(2X_{2i}) + 1.4(\exp(-0.2X_{1i}) + \exp(0.7X_{2i}^2)) + \varepsilon_i$.

Bandwidths and kernel function used in the simulation were taken to be $h = n^{-1/4.5}, w = n^{-1/6}(\log n)^{1/6}$ and

$$K(z) = K^*(z) = \tfrac{9}{16}(1 - z_1^2)(1 - z_2^2)I(|z_1| \leq 1, |z_2| \leq 1).$$

The sample sizes chosen are $50, 100, 200$ and $300$, each repeated 1000 times. Table 3 lists means and MSE of $\hat{\theta}_n = (\hat{\theta}_{n1}, \hat{\theta}_{n2})'$ obtained by minimizing $M_n(\theta)$ and employing the Newton–Raphson algorithm. As in case 1, one sees little bias in the estimator for all chosen sample sizes.

Table 4 gives the empirical sizes and powers of the $\widehat{\mathcal{D}}_n$-test for testing Model 0 against Models 1–3. From this table one sees that this test is conservative when sample sizes are small, while empirical levels increase with the sample sizes and indeed preserve the nominal size 0.05. It also shows that the MD test performs well for sample sizes 200 and larger at all alternatives.

TABLE 3
Mean and MSE of $\hat{\theta}_n$

| Sample size | 50 | 100 | 200 | 300 |
|---|---|---|---|---|
| Mean of $\hat{\theta}_{n1}$ | 0.9978 | 0.9973 | 0.9974 | 0.9988 |
| MSE of $\hat{\theta}_{n1}$ | 0.0190 | 0.0095 | 0.0053 | 0.0034 |
| Mean of $\hat{\theta}_{n2}$ | 1.9962 | 1.9965 | 2.0013 | 2.0004 |
| MSE of $\hat{\theta}_{n2}$ | 0.0063 | 0.0028 | 0.0014 | 0.0010 |

TABLE 4
Levels and powers of the minimum distance test in case 2

| Sample size | 50 | 100 | 200 | 300 |
|---|---|---|---|---|
| Model 0 | 0.003 | 0.019 | 0.049 | 0.052 |
| Model 1 | 0.158 | 0.843 | 0.979 | 0.996 |
| Model 2 | 0.165 | 0.840 | 0.976 | 0.992 |
| Model 3 | 0.044 | 0.608 | 0.954 | 0.997 |

**Acknowledgment.** Authors would like to thank the two referees and Jianqing Fan for constructive comments.

DEPARTMENT OF STATISTICS
  AND PROBABILITY
MICHIGAN STATE UNIVERSITY
EAST LANSING, MICHIGAN 48824-1027
USA
E-MAIL: koul@stt.msu.edu

DEPARTMENT OF STATISTICS
KANSAS STATE UNIVERSITY
MANHATTAN, KANSAS 66506-0802
USA
E-MAIL: weixing@ksu.edu